\documentclass[12pt]{article}
\newfont{\inhead}{eufm10}

\newcommand{\pbar}{{\bar\partial}}

\newcommand{\Z}{{\mathbb Z}}
\newcommand{\R}{{\Bbb R}}
\newcommand{\C}{{\mathbb C}}

\newcommand{\Q}{{\Bbb Q}}

\newcommand{\Aa}{{\cal A}}   

\newcommand{\Ee}{{\cal E}}
\newcommand{\Ff}{{\cal F}}
\newcommand{\Gg}{{\cal G}}   

\newcommand{\Mm}{{\cal M}}   

\newcommand{\Ss}{{\cal S}}
\newcommand{\Kk}{{\cal K}}

\newcommand{\coker}{{\rm coker }}  
\newcommand{\im}{{\rm im }}        
\newcommand{\diag}{{\rm diag}}     

\newcommand{\Diff}{{\rm Diff}}        
\newcommand{\Symp}{{\rm Symp}}        
\newcommand{\Map}{{\rm Map}}          
\newcommand{\End}{{\rm End}}          

\newcommand{\D}{{\rm D}}

\newcommand{\point}{{\rm pt}}

\newcommand{\G}{{\rm G}}

\newcommand{\U}{{\rm U}}

\newcommand{\PSL}{{\rm PSL}}

\newcommand{\Cinf}{C^{\infty}}

\def\NABLA#1{{\mathop{\nabla\kern-.5ex\lower1ex\hbox{$#1$}}}}
\def\Nabla#1{\nabla\kern-.5ex{}_#1}

\newcommand{\p}{{\partial}}

\newcommand{\proof}[1]{\noindent{\bf Proof#1:\  }}
\newcommand{\obs}[1]{\noindent{\bf Observation#1:\  }}

\newcommand{\QED}{\hfill$\square$\medskip}

\newtheorem{theorem}{Theorem}[section]
\newtheorem{corollary}[theorem]{Corollary}

\newtheorem{lemma}[theorem]{Lemma}
\newtheorem{proposition}[theorem]{Proposition}
\newtheorem{definition}[theorem]{Definition}

\usepackage{amssymb}
\usepackage{amsmath}
\usepackage[all]{xy}

\begin{document}

\title{Relative family Gromov-Witten invariants  and symplectomorphisms}
\author{Olgu\c ta Bu\c se\footnote{SUNY at Stony Brook, {\tt buse@math.sunysb.edu}}}
\date{\today}
\maketitle

\abstract{\it
 
We study the symplectomorphism groups $G_{\lambda}=\Symp_0(M,\omega_{\lambda})$ of an arbitrary closed manifold M equipped with a 1-parameter family of symplectic forms $\omega_{\lambda}$ with variable cohomology class. We show that the existence of nontrivial elements in $\pi_*(\Aa,\Aa')$, where $(\Aa,\Aa')$ is a suitable pair of spaces of almost complex structures, implies the existence of  families of nontrivial elements in  $\pi_{*-i}G_{\lambda}$, for $i=1$ or $2$. Suitable parametric Gromov Witten invariants detect nontrivial elements in $\pi_*(\Aa,\Aa')$. By looking at certain resolutions of quotient singularities we investigate the situation $(M,\omega_{\lambda})=(S^2 \times S^2 \times X,\sigma_F \oplus \lambda \sigma_B \oplus \omega_{st})$, with $(X,\omega_{st})$ an arbitrary symplectic manifold. We find families of nontrivial elements in  $\pi_k(G_{\lambda}^X)$, for  countably many $k$ and different values of $\lambda$. In particular we show that the fragile elements $w_{\ell}$ found by Abreu-McDuff \cite{AM} in $\pi_{4 \ell}(G_{\ell+1}^{\point})$ do not disappear when we consider them in $S^2 \times S^2 \times X$.
}

\vspace{1cm}
\section{Introduction}
  Consider $(M^{2n},\omega)$ a $2n$ dimensional compact symplectic manifold.
A basic invariant which distinguishes among different symplectic structures 
on $M$ is the group of symplectomorphisms, $\Symp(M,\omega)$. This is an infinite dimensional group  endowed with a natural $C^{\infty} $ topology.

Two natural questions arise in relation with  $\Symp(M,\omega)$  namely
\begin{description}
\item[(1)] What can be said about the topological type of  $\Symp(M,\omega)$?

\item[(2)]
 How does the topological type change as $\omega$ varies?
\end{description}

Research has been done in this direction by various authors [Abreu \cite{A}, Le-Ono \cite{LO}, McDuff \cite{M3}, Seidel \cite{SE}] by using information on $J$-holomorphic curves.
We investigate these questions by defining relative parametric
GW invariants, which are sensitive to the topology of appropriate spaces of 
almost complex structures. The connection between the spaces of almost 
complex structures and the symplectomorphism groups is achieved by means of
the following fibration, introduced by Kronheimer \cite{K} and used in McDuff \cite{M4}:

\begin{equation}\label{fr:Kseq}  
\xymatrix
{
  \Symp_0(M,\omega) \ar[r] & 
    \Diff_0(M) \ar[rrr]^{\psi \rightarrow (\psi^{-1})^*\omega} & & & S_{\omega}
}
\end{equation}
where $S_{\omega}$ is  the space of symplectic forms which can be joined to
$\omega$ through a path of cohomologous symplectic forms,  $ \Diff_0(M)$ is
the connected component of the identity inside the group of diffeomorphism 
and $\Symp_0(M,\omega) = \Symp(M,\omega)\cap \Diff_0(M)$. Now consider $\Aa_{\omega}$  the space of all almost complex structures tamed by some 
symplectic form $\omega'$ in   $\Ss_{\omega}$. By \cite{M4},   $\Aa_{\omega}$ is homotopy equivalent to  $\Ss_{\omega}$. This yields  the following 
homotopy fibration:

\begin{equation}\label{fr:Dseq} 
\xymatrix
{
  \Symp_0(M,\omega) \ar[r] & \Diff_0(M) \ar[r] & \Aa_{\omega}.
 }
 \end{equation}
Our strategy will be to define suitable pairs $(\Aa,\Aa')$ of spaces of almost complex structures, such that information on nontrivial homotopy groups in $(\Aa,\Aa')$ extends to  information on  $\Symp_0(M,\omega)$. We develop a version of relative GW invariants in family which detects such nontrivial elements in $\pi_*(\Aa,\Aa')$.

\noindent {\bf Outline of the methods}

In section {\bf 2} we  will define the invariants as follows:
Consider $D \in H_2(M,\Z)$ and let $\Aa_{\omega,\D }^c$ be the subspace of 
$\Aa_{\omega}$ consisting of those almost complex structures $J$ which  do not  admit 
$J-$holomorphic stable maps in the class $D $. For $I$ an interval in $\R$ further define $(\Aa_I,\Aa_{I,D}^c)=\bigcup_{\omega_{\lambda} \in L}(\Aa_{\omega_{\lambda}},\Aa_{\omega_{\lambda},D }^c)$, where the cohomology of symplectic form  $[\omega_{\lambda}] $ is deformed along a line L inside a positive cone $\Kk \in H^2(M,\R)$.
 
Consider  a family of almost complex structures $(J_B,J_{\p B})$  that represent an element in $\pi_*(\Aa_{I},\Aa_{I,D}^c).$
We will define a homomorphism 
\begin{equation}
PGW^{M,(J_B,J_{\p B})}_{D,0,k}:\bigoplus_{i=1}^k H^{a_i}(M,\Q)^k \rightarrow \Q
\end{equation}
by counting $J_b$-holomorphic stable maps in class D, for all $b \in B$. This is well defined because the class D is never represented as a $J_b$-holomorphic stable maps if $b \in \p B$.
We have the following
\begin{theorem}
{\bf i)} The  invariants $ PGW^{M,(J_B,J_{\p B})}_{D,0,k}$ are symplectic deformation invariants and  depend only on the relative homotopy class of the pair $(J_B,J_{\p B})$.

{\bf ii)} For  a fixed choice of $k,D$ and $\alpha_i$ the map $\Theta _{0,k,\alpha_1,\ldots,\alpha_k}:\pi_*(\Aa_I,\Aa_{I,D}^c) \rightarrow \Q$, given by
$$\Theta _{k,\alpha_1,\ldots,\alpha_k}([(J_B,J_{\p B})] = PGW^{M,(J_B,J_{\p B})}_{D,0,k}(\alpha_1,\ldots,\alpha_k).$$
 is a homomorphism.
\end{theorem}
The reason why {\bf (i)} holds is that the  class D is never represented for a $J_b$ with  $b \in \p B$.

 In section {\bf 3}  we will exhibit some  examples of  nontrivial PGW. There we consider the case when  $ (M,\omega)$ is $S^2 \times S^2 \times X $, where X is an arbitrary symplectic manifold and $\omega =\omega_{\lambda} \oplus \omega_{st}$, with 
$\omega_{\lambda}=\sigma_F \oplus \lambda \sigma_B$. Here
$\sigma_F,\sigma_B$ are forms on the  fiber and  base respectively, of total area 1,$\lambda \geq 1$,  and $\omega_{st}$ is arbitrary symplectic form on X.
The families $(J_B,J_{\p B})$  of almost complex structures  are provided for $S^2 \times S^2$ in \cite{K} and then further investigated in \cite{AM}. One has to look at a
quotient singularity, $\C^2/C_{2\ell }$, where $C_{2\ell }$ is the cyclic group 
of order $2\ell $ acting diagonally  by scalars on $\C^2.$ The deformation space for the canonical resolution of this singularity provides a $4\ell -2$ family 
$(J_{B_{\ell}},\p J_{ B_{\ell}}) \in  (\Aa_{[\ell+\epsilon,\ell]},\Aa_{\ell})$ for 
which suitable PGW are nontrivial.

The link between these examples and the corresponding groups of symplectomorphisms will be explained in section {\bf 4}. It will be there where we explain the extent to which the known homotopy properties (see
\cite{AM})  of  $\Symp_0(S^2 \times S^2, \omega_{\lambda})$ are reflected 
in the high homotopy groups of  $ G_{\lambda}^X :=\Symp_0(S^2 \times S^2 \times X,\omega_{\lambda} \oplus \omega_{st})$. For every $M,\omega_{\lambda}$ a general symplectic manifold, we set the notation $G_{\lambda}:=\Symp_0(M,\omega_{\lambda})$.\\
 To be able to give any answers related to the two questions posed in the beginning, one has to establish first a more precise language in which they make sense. One of the difficulties is that in general there is no direct map  $G_{\lambda} \rightarrow G_{\lambda + \epsilon}$. In the particular situation $M=S^2 \times S^2 \times \point$  Abreu-McDuff in \cite{AM} and \cite{M4} find natural maps  $G_{\lambda}^{\point}\rightarrow G_{\lambda + \epsilon}^{\point}$, well defined up to homotopy, and prove:
\begin{theorem}{(Abreu-McDuff)}\label{th:mcduffabreu}
(i) The homotopy type of $G_{\lambda}^{\point}$ is constant on all the intervals $(\ell-1,\ell]$ with $\ell \geq 2 $ a natural number. Moreover, as $\lambda$ passes an integer $\ell,\ell \geq 2$ the groups $\pi _i (G_{\lambda}^{\point}),  i\leq 4 \ell -5$,  do not change. 

(ii) There is an element  $w _{\ell} \in \pi_{4 \ell-4} (G_{\lambda}^{\point}) \times \Q $ when $\ell-1 < \lambda \leq \ell $ that vanishes for $\lambda > \ell $.
\end{theorem}

When we deal with a general manifold M, to get around the fact that there is no map $G_{\lambda} \rightarrow G_{\lambda + \epsilon}$ we show that for any compact $K \subset G_{\lambda}$, the inclusion $0 \times K \subset G_{\lambda}$ extends  to a map $h$ that fits into the following commuting diagram:
\begin{equation}
\xymatrix{
  h:[ - \epsilon, \epsilon] \times K  \ar[rr] \ar[d]^{pr_1} &  &  \Gg:=\bigcup(G_{\lambda} \times \lambda) \subset \Diff \times \R \ar[d]^{pr_2} \\
  [ - \epsilon, \epsilon] \;\;\ar[rr]^{incl}  & &\;\;\;\;\;(- \infty, \infty) }
\end{equation}
Moreover, for any two  such  maps $h$ and $h'$ which coincide on $0 \times K$, there is, for $\epsilon ' $ small enough, a homotopy between them $H:[0,1] \times [- \epsilon ',\epsilon '] \times K \rightarrow \Gg$ which also preserves the fibers of the natural projections.
We therefore  see   that, for any  cycle $\rho$ in  $G_{\lambda}$ there are {\it extensions} $\rho_{\epsilon}$ in  $G_{\lambda + \epsilon}$ which, for  $\epsilon$ sufficiently small, are unique  up to homotopy. Hence they give well defined elements in $\pi_*G_{\lambda + \epsilon}$. 
 
It will therefore make sense to ask what will become of an element $\rho \in \pi_*G_{\lambda}$ inside  $\pi_*G_{\lambda + \epsilon}$, for small $\epsilon$. In this language we say that an element $\theta_{\ell} \in \pi_*G_{\ell}$ is {\it fragile} if any  extension $\theta_{\ell+\epsilon}$ is null-homotopic  in $\pi_*(G_{\ell + \epsilon})$ for $\epsilon>0$. Also, we say that a family $\eta_{\ell+\epsilon} \in \pi_*G_{\ell+\epsilon},0<\epsilon$ is {\it new} if there is no $\eta_{\ell} \in \pi_*G_{\ell}$ whose extension is $\eta_{\ell+\epsilon}$. We consider the space $\Aa_{\ell^+}$ roughly given by $\Aa_{\ell^+}:=\bigcap_{0<\epsilon<\epsilon_0}\Aa_{\ell+\epsilon}$.\\ We say that an element $\alpha \in \pi_*(\Aa_{\ell^+},\Aa_{\ell})$ is {\it persistent} if it has nonzero image under the map
 $\pi_*(\Aa_{\ell^+},\Aa_{\ell}) \rightarrow \pi_*(\Aa_{[\ell,\ell+\epsilon]},\Aa_{\ell})$.\\
The content of our main theorem is the following:

\begin{theorem}\label{th:Mts}
 Assume that we have a  persistent element  $0 \ne \beta_{\ell} \in \pi_k (\Aa_{\ell^+},\Aa_{\ell})$ 
  Exactly one  of the
statements below holds.

{\bf A)} There is an uniquely associated  non-zero {\bf fragile} element $\theta_{\ell} \in \pi_{k-2}G_{\ell}$, 
  such that $i_*(\theta_{\ell})=0 $ in $\pi_{k-2}\Diff_0 (M )$. 

{\bf B)} There exists an $\epsilon_{\ell} >0$ and an uniquely associated family of {\bf new} elements $0 \ne \eta_{\ell +\epsilon} \in \pi_{k-1}G_{\ell+\epsilon},0<\epsilon<\epsilon_{\ell}$.
\end{theorem}

 We should point out that our methods do not allow us to decide in general whether the image of $\eta_{\ell+\epsilon}$  in $\Diff_0(M)$ is zero or not. 

 We show that the hypothesis of the theorem is verified  when $M=S^2 \times S^2 \times X$. We consider $D=A - \ell F$. Since $[\sigma_F \oplus \lambda \sigma_B \oplus \omega_{st}](A-\ell F)=0$ we get that $\Aa_{\ell} \subset \Aa_{[\ell,\ell+\epsilon],D}^c$. In this situation the $4 \ell -2$ dimensional elements $(B_{\ell},\p B_{\ell})$ obtained in section {\bf 3}  are  detected as nontrivial in $\pi_{4\ell-2}(\Aa_{\ell^+},\Aa_{\ell})$ and are persistent. In fact in general PGW invariants detect persistent elements.
 By varying the value of the integer $\ell$ we obtain infinitely many  values of $\lambda$ for which higher order homotopy groups of $G_{\lambda}^X$ are  nontrivial and also make a more detailed
 discussion regarding the stability of the elements $w_{\ell}$ provided by theorem \ref{th:mcduffabreu} inside $G_{\lambda}^X$. This is the content of the following:
  \begin{corollary}\label{cr:Mts}
  For any natural number $\ell \geq 1$, exactly one of the
statements below holds.

{\bf A)} There is a  non-zero {\bf fragile} element $w_{\ell}^X \in \pi_{4\ell-4}G_{\ell}^X$, 
  such that $i_*(w_{\ell}^X)=0 $ in $\pi_*\Diff_0 (S^2 \times S^2 \times X )$. This element can be identified with $w_{\ell} \times id$.

{\bf B)} There exists an $\epsilon_{\ell} >0$ and a family of {\bf new} elements $0 \neq \eta_{\ell+\epsilon}^X \in \pi_{4\ell-3}G_{\ell+\epsilon}^X,0<\epsilon<\epsilon_{\ell}.$ 
\end{corollary}

In particular this shows that the fragile elements obtained by Abreu-McDuff for$\ell>1$ do not disappear when we consider them inside $S^2 \times S^2 \times X$.  One possibility is that $0 \neq w_{\ell}\times id \in \pi_{4 \ell-4} (G_{\ell }^{X})$. If this is not the case, then we have the associated new $4 \ell-3$ dimensional elements $0 \neq \eta_{\ell+\epsilon}^X$ in $\pi_{4\ell-3}G_{\ell + \epsilon}^{X}$ for small $\epsilon>0$.
 For general X, when $\ell=1$ it is known by work of Le-Ono that {\bf B} takes place and  $0 \neq i_*(\eta_{\ell+\epsilon}) \in \Diff(S^2 \times S^2 \times X)$. Also, for $X=\point$ and $\ell>1$ from the work of Abreu-McDuff we know that {\bf A} takes place. 

We do not have examples when case {\bf B} takes place and $i_*(\eta_{\ell+\epsilon}) \neq 0 \in \Diff(M).$

Our method had been inspired by P. Kronheimer's work, who uses parametric 
Seiberg-Witten invariants in dimension 4, as well 
as by the work of D. McDuff \cite{M3}.  Similar work has been done in this 
direction by Le-Ono in\cite{LO}; by looking at related but slightly different  parametric GW 
invariants they get results about $\pi_i(\Symp_0(S^2 \times S^2 \times X,\omega_{\omega_1 \oplus \omega_{st}}))$ when $i=1,3.$
 In section {\bf 3} we could consider $\C^2 / \C_{2 \ell+1}$ instead and by carrying out similar arguments get the same type of results for $\C P^2 \#\overline{\C P^2} \times X$. 

\noindent{\bf Acknowledgments} This is part of the author's doctoral research at SUNY Stony Brook. The author would like to thank her advisor, Dusa McDuff, for her suggestions, advice and comments on earlier drafts.

\section{Relative parametric GW invariants }\label{sec2}
\subsection{Definition and properties}\label{subsec21}

Consider B to be a compact manifold with boundary  and a  smooth map $i:(B,\p B) \rightarrow   (\Aa_I, \Aa_{I,D}^c)$. Although the invariants can be defined with this data regarding the parameter space,  for the applications we have in mind we will consider $B$ to be an n-ball such that $i$ represents a relative homotopy class
 in  $\pi_*(\Aa_ I, \Aa_{I,D}^c)$. We will often write  $J_b:=i(b)$ and $J_B=\im(i)$, 
and refer to $\im B$ in $\Aa_I$ as $J_B$. Consider also a smooth family of symplectic forms $(\omega_b)_{b \in B}$ where $\omega_b$  tames $J_b$. We point out that the
$\omega_b$ need not be cohomologous, as the taming condition is an open condition. Our goal here 
is to show how  we can define parametric GW invariants relative to the boundary $\p J_B$, which count $J_b$ holomorphic maps for some $b \in B$. These will not depend  either on deformations of the family $\omega_B$ or on the relative homotopy class   $(J_B,\partial J_ B) \subset 
   (\Aa_{I},\Aa_{I,D}^c).$

 Consider $\widetilde \Mm_{0,k}^*(M,D,(J_B,\p J_B))$ the space of tuples $(b,f,x_1,\ldots, x_k)$ where $f:S^2 \rightarrow M$ is a simple \footnote{We say that $f:\Sigma \rightarrow M$ is simple if it is not the composite of a holomorphic branched covering map $(\Sigma, j) \rightarrow (\Sigma',j') $ of degree greater than 1 with a J-holomorphic map $\Sigma' \rightarrow M$.}    $J_b$-holomorphic map in class D, for some $b \in B$ and $x_i$ are pairwise distinct points on $S^2$.          
We will consider $$ \Mm_{0,k}^*(M,D,(J_B,\p J_B))= \widetilde \Mm_{0,k}^*(M,D,(J_B,\p J_B))/G$$
 where $G=\PSL(2,\C)$ acts on the moduli space by reparametrizations of the domain. Denote the elements of $\Mm_{0,k}^*(M,D,(J_B,\p J_B)$ by  $[b,f,x_1,\ldots,x_k]$.

In the best scenario, for a good choice of $(J_B,\p J_B)$, the following hold:\\
 {\bf (1)} $\widetilde \Mm_{0,k}^*(M,D,(J_B,\p J_B))$ is a manifold of dimension $2n+2c_1(D)+2k+\dim B$ and \\
{\bf (2)} $ \Mm_{0,k}^*:=\Mm_{0,k}^*(M,D,(J_B,\p J_B))$ is compact. \\
Then the image of the map 
\begin{equation}\label{eval}
ev:\Mm_{0,k}^*(M,D,(J_B,\p J_B)) \rightarrow M^k
\end{equation}
with $ev([b,f,x_1,\ldots,x_k]):=(f(x_1),\ldots,f(x_k))$ will provide a cycle $ ev_*(\Mm_{0,k}^*)$ in $M^k$ which, by intersection with homology classes of complementary dimension in $M^k$, gives the parametric Gromov-Witten invariants.

As we will see in the regularity discussion below, {\bf (1)} is always possible to accomplish by Sard-Smale theorem. However, {\bf (2)} is seldom true; the compactification $ \overline \Mm_{0,k}(M,D,(J_B,\p J_B))$  of $ \Mm_{0,k}^*(M,D,(J_B,\p J_B))$ contains both {\it stable maps} and {\it nonsimple curves}, which we sometimes call multiple cover curves. We will spell out some of these notions later in this section; for more information, the reader can check \cite{MS}, \cite{LT}, \cite{R}, \cite{CO}.

  $ \overline \Mm_{0,k}(M,D,(J_B,\p J_B))$ is a stratified space. The best we can hope is that the image of the
 evaluation map $ev:\Mm_{0,k}^* \rightarrow M^k$  a {\it pseudo-cycle}, or differently said, the boundary strata in the image will be codimension 2 or bigger. 
If this scenario works we can still define the PGW as the intersection between the image of $ev$ and classes of complementary dimension in $H_*(M^k)$. 
This will be for instance the case when the class $D$ is $J_b$ indecomposable for any $b \in B$, that is, no $J_b$ holomorphic map in class D can decompose into a connected union of $J_b$ holomorphic spheres $C=C^1 \bigcup C^2 \bigcup \ldots \bigcup C^N$ such that each $C^i$ represents the class $D_i$ and $D=D_1+ \ldots +D_N$.  In fact, in this situation the image of $ev$ is a cycle. This hypothesis will be enough for the application we have in mind. 

\begin{definition}
We will say that the {\bf hypothesis  $H_1$ is satisfied} if the class $D$ is $J_b$ indecomposable  for any $b \in B$. 
\end{definition}

\noindent{\bf Parametric regularity}

We begin by explaining what is {\it D-parametric regularity} and contrast it with the usual D-regularity for $J$ (see \cite{MS}). For this we need to introduce the following facts.

Let ${\cal X}=\Map(\Sigma,M;D)$ be the  space of {\it somewhere injective} \footnote{We say that a map $  f:\sigma \rightarrow M$ is somewhere injective if $df(z) \neq 0, f^{-1}(f(z)={z}$ for some $z \in \Sigma$. A simple J-holomorphic map is somewhere injective (see \cite{MS}.} smooth maps  $f:\Sigma \rightarrow M$ representing class D.
 This is an infinite dimensional manifold with $T_f {\cal X} = \Cinf (f^*TM)$. 
We will next  consider the following generalized vector bundle
${\cal E} \longrightarrow B \times {\cal X}$, whose fiber at $(b,f)$ is the 
space ${\cal E}_{b,f}=\Omega_{J_b}^{0,1}(f^*TM)$ of smooth $J_b$ antilinear 
forms with values in $f^*TM$. In this vector bundle we  consider a section $\Phi :B \times {\cal X} \longrightarrow {\cal E}$, given by 
\begin{equation}\label{phi}
\Phi(b,f)=\frac12(df+J_b \circ df \circ j)
\end{equation}

The zeros of $\Phi$ are precisely $J_b$ holomorphic maps and
thus the moduli space  
 $$\widetilde \Mm_{0,0}^*(M,D,(J_B,\p J_B))=\Phi^{-1}(0),$$
is  the intersection of $\im \Phi$ with the zero section of the bundle. 
Since we would like \\
 $\widetilde \Mm_{0,k}^*(M,D,(J_B,\p J_B))$ to be a manifold we require that 
$\Phi$ is transversal to the zero section. This means that the image of 
$d\Phi(b,f)$ is complementary to the tangent space $T_b B\oplus T_f \cal {X}$ 
of the zero section. 
But for any $f$ which is $J_b$ holomorphic, $d \Phi$ is given by
$$d \Phi(b,f):T_b B \oplus \Cinf (f^*TM) \longrightarrow 
  T_b B  \oplus T_f{ \cal X} \oplus {\cal E}_{b,f}$$
 If we consider now the projection onto the vertical space of the bundle:
$$proj_2: T_bB \oplus T_f {\cal X} \oplus {\cal E}_{b,f} \longrightarrow 
  {\cal E}_{b,f}$$
the above transversality translates into the fact that
\begin{equation}\label{dphi}
d \Phi(b,f) \circ proj_2:T_bB \oplus \Cinf (f^*TM) \longrightarrow
  \Omega_{J_b}^{0,1}(\Sigma,f^*TM)
\end{equation}
is onto. We will make the notation $D \Phi(b,f)=d \Phi(b,f) \circ proj_2$. We 
then have:
\begin{definition}
We say that a $J_b$ holomorphic map f is $J_B$ parametric 
regular if $D \Phi (b,f)$ is onto.
\end{definition}

\obs{} 
The linearized operator is well defined if there is no pair $(b,f)$ with 
$f$  $J_b$ holomorphic and $b \in \partial B$. This is precisely the condition 
we imposed on  $(J_B,\p J_B)$ to give  a relative cycle  in
$(\Aa_{I},\Aa_{I,D}^c).$

\begin{definition}
Consider $(J_B,\omega_B)$ as above. We say that  $(J_B,J_{\p B})$ is an D-parametric regular family 
of almost complex structures if any $J_b$  holomorphic map in class D is 
 parametric regular. We denote by $J_{preg}(D)$ the set of all D-parametric regular families $(J_B,\partial J_ B) \in  (\Aa_I,\Aa_{I,D}^c)$.
\end{definition}

In order to apply the implicit function theorem and Sard-Smale theorem we must work on Banach manifolds and hence complete all spaces under suitable Sobolev norms. For example, one should to work on  spaces  consisting of almost complex structures of class $C^l$, on ${\cal X}^{k,p}$, with $kp>2$, the space of maps whose k-th derivatives are of class $L^p$. Also, we should work on  
$${\cal E}_f^p=L^p(\Lambda ^{0,1} \otimes _J f^* TM))$$ rather that with $\Omega^{0,1}_J(\Sigma,f^*TM)$.
 There are standard arguments \cite{MS} to show that one can  pass the following arguments from spaces of $C^l$ objects (which are Banach manifolds)  to spaces of $C^{\infty}$ objects  (which are Frechet manifolds). For simplicity we will drop the superscripts $l,k,p$ unless it  will be relevant to specify them. 
We have the following:

\begin{theorem}\label{th:regexs}

 If $J_B \in J_{preg} (D)$, then the moduli space $\widetilde \Mm_{0,0}^*(M,D,(J_B,\p J_B))$ is 
a smooth open manifold of dimension $2n + 2 c_1(D) + \dim B$, with a natural 
orientation.
\end{theorem}

Moreover,  if one considers $\widetilde \Mm_{0,0}^*(M,D,(J_B,\p J_B)) \times (S^2)^k $ and takes away all the diagonals of the type
 $\widetilde \Mm_{0,0}^*(M,D,(J_B,\p J_B)) \times \diag_{i,j}$, what we obtain is precisely $\widetilde \Mm_{0,k}^*(M,D,(J_B,\p J_B))$. This will therefore be a manifold of dimension $2n +2c_1(D)+\dim B +2k$.

Let  $\widetilde \Mm_{0,0}^*(M,D,\Aa_I)$ be the {\it universal} moduli space consisting of pairs $(f,J)$ where $J \in \Aa_I$ and $f$ is $J$ -holomorphic. It will be more relevant to the story to point out the following characterization of parametric regularity.
\begin{proposition}\label{pr:transreg}
Consider the diagram
\begin{equation}
\xymatrix{
   & \widetilde \Mm_{0,0}^{ *}(M,D,\Aa_I)  \ar[d]^{ \Pi } \\
  (B, \partial B) \ar[r]^i & ( \Aa_I,\Aa_{ I, D}^{c}) }
\end{equation}
Then $J_B \in J _{preg} (A) $ iff $ i \pitchfork \Pi. $ 
\end{proposition}
\proof{} For simplicity we will denote by $D_{f,b}=D \Phi(b,f)_ {| C^{ \infty} (f^*(TM)}$. By (\ref{dphi}) the surjectivity of $D \Phi(b,f)$ is then equivalent with the surjectivity of the following linear operator
$$D \phi _{| T_b B}:T_b B \rightarrow \coker D_{b,f}$$
We will denote $i(b)=J$. The tangent space $T_J \Aa_{I}$ to $\Aa_{I}$ consists of all  sections Y of the bundle $\End(TM,J)$ whose fiber at $p \in M$ is the space of linear maps $Y :T_p M \rightarrow T_p M$ such that
$YJ+JY=0$; we will consider the map
$$R :T_J \Aa_{I} \rightarrow \Omega^{0,1}_J \left(\Sigma,f^*TM \right)$$
given by $R(Y)=\frac 12 Y \circ df \circ j.$ The map
 $$d \Pi:T_{f,J} \widetilde \Mm_{0,0}^{*}(M,D,\Aa_I)
\rightarrow T_J \Aa_{I}$$
is given by $d \Pi(\xi,Y)=Y$, where the pair $(\xi,Y)$  is in $T_{f,J} \widetilde \Mm_{0,0}^{*}(M,D,\Aa_I)$ if and only if 
\begin{equation}\label{defpi}
 D_{f,b}(\xi)+R(Y)=0
\end{equation}
From this one can see that  $\im D_{f,b} = R(\im (d\Pi))$. Since $D_{b,f}$ is elliptic and $\ker R \subset \im d \Pi$, it follows that $\coker d\Pi $ has finite dimension. If we  consider the map ${\cal F}:{\cal X} \times \Aa_I \rightarrow {\cal E}$, given by ${\cal F}(f,J)=\pbar_J(f)$ then (see \cite{MS}) the linearization  at a zero $(f,J)$ with $f$ simple is onto. That is 
$$D{\cal F}(f,J)(\xi,Y)=D_f \xi + R(Y)$$ is onto.This implies  that  
$\coker D_f$ is covered by $R$. We can show that there is an induced map
 $$\widetilde{R}:\coker d \Pi \rightarrow \coker D_{b,f}$$ which is isomorphism. The proof of the proposition then follows easily. 
$D \Phi _{| T_b B}(Y)=R \circ di$, so we have 
 $i \pitchfork \Pi \Leftrightarrow di \rightarrow 
  \coker d\Pi \;\text{ onto }  \Leftrightarrow 
  \widetilde{R} \circ di \rightarrow \coker D_{b,f} \; \text{ onto.}$
\QED

 There are few key points to be noticed here. Notice that parametric regularity is a generalization of the usual regularity. Indeed, if we consider $J_b=J$ to be constant for $b$ in a neighborhood around $b_0$  then the regularity of an almost complex structure $J$  simply says, following the diagram above,
that $d \Pi$ is surjective. If we now regard $J$  within an arbitrary family $J_B$, this no longer needs to be the case. It will then suffice that the cokernel of $d \Pi^l$ is covered by the variation of $J$  in the direction of B. In subsection \ref{subsec22} we will see that, when we count rational maps, an equivalent criterion of parametric regularity will be given by the usual regularity in some suitable ambient space.

 There exist a large  subset of parametric regular families of almost complex structures inside $(\Aa_I,\Aa_{I,D}^c)$. This is because one can employ Sard-Smale theorem \cite{S} and show that any map $i:(B, \p B) \rightarrow (\Aa_I,\Aa_{I,D}^c)$ in the prop (\ref{pr:transreg}) can be perturbed to an $i'$ such that
 $i' \pitchfork \Pi$.  
\begin{definition}
We will say that $(J_B,\p J_B)$ {\bf satisfies hypothesis  $H_2$} if it is a {\it D-parametric regular} family of almost complex structures.
\end{definition}

\noindent{\bf Compactness}

 We have already mentioned that  {\bf $H_1$ }  is not verified for all parametric regular families  of  almost complex structures. Moreover, we cannot guarantee that for an arbitrary choice of a regular family, the image of the evaluation map (\ref{eval}) is a pseudocycle. That is because we
 could potentially have nonsimple elements in the compactification $ \overline \Mm_{0,k}(M,D,(J_B,\p J_B)$ which would yield boundary strata of high 
dimension.

In the situation that $B=\point$ there are various procedures [Li-Tian (\cite{LT}), Ruan (\cite{R}), Fukaya-Ono (\cite{FO})  to build up a theory which would provide {\it a virtual moduli cycle}, that is, an object which carries a fundamental class required for the definition of the invariants.

Roughly speaking, locally one needs to consider here all the stable holomorphic  maps as well as small perturbations of these. There are then various
 procedures to pass to a global object with the required properties.
These go through without essential changes if one considers parameter spaces with no boundary [see Leung-Bryan(\cite{CO}, Ruan(\cite{R})].

To make this more precise, let's first give the following 
\begin{definition}\cite{LT}
 A stable smooth  rational map is given by a tuple 
$(f,\Sigma,x_1,\ldots,x_k)$ satisfying:

1) $\Sigma=\bigcup_{i=1}^m \Sigma_i$ is a connected rational curve with normal crossing singularities  and $x_1,\ldots,x_k$ are distinct smooth points in $\Sigma$

2) f is continuous and each restriction $ f_{| \Sigma_{i}}$ lifts to a smooth map from the normalization $\overline{\Sigma_{i}}$   to M;

3) If $f_{|\Sigma_{i}}$ is constant then $\Sigma_i$ contains at least three  $\it{special}$ points. Here, a special point is either a singular point or a marked point.
\end{definition}

We will quotient the space of stable maps by the group of automorphisms of the domain; define $\overline{\Ff}_D(M,0,k)$ to be the space of equivalence classes $ [f,\Sigma,x_1,\ldots,x_k]$.
The bundle ${\cal E}$ can be then viewed as a generalized bundle over $\overline{\Ff}_D(M,0,k)$ and together with the section $\Phi$ defines a generalized Fredholm bundle in the sense of Li-Tian, of index $2n +2c_1(D)+2k -6$.

In  our situation we need to make sure that the boundary causes no problem.
The following  lemma basically states that if we consider an  appropriately small open  neighborhood  of  $ \overline \Mm_{0,k}(M,D,(J_B,\p J_B)$ consisting
 of almost  holomorphic maps, then  its projection onto  $J_B$ stays away from $\p J_B$.
\begin{lemma}
For any compact set $J_B \in \Aa_I $ such that $\p J_B \subset  \Aa_{I,D}^c$  $  \exists $ a $\delta >0$ and $\epsilon (\delta) >0  $ for which  there is no stable map  $[f,\Sigma,x_1,\ldots,x_k]$ such that  $\pbar_J f =\nu$, when  $d(J,  \p J_B) < \delta$ and  $\nu \in L^p(\Lambda ^{0,1} \otimes _J f^* TM)$  with $ |\nu| \leq \epsilon (\delta) $. 
\end{lemma}

\proof{} We will prove this by assuming the opposite.
Assume that we have a sequence $J_i$,  $\nu _i$ and $ f_i$ such that $d(J_i,\p J_B) \rightarrow 0$,  $| \nu _i | = \epsilon _i  \rightarrow 0$ and each $f_i$ is a stable map in class D with the property that   $\pbar_{J_i} f_i =\nu _i$.
Since $J_B$ is compact we find a convergent subsequence $J_i$, 
whose limit $J_{\infty}$  is in $\p J_B$. But this would lead to a
 contradiction because by the Gromov compactness theorem there is a subsequence  of $f_i$ which converges to a $J_{\infty}$
 stable holomorphic map in class D. This will contradict the fact that  $  J_{\infty} \in \p J_B \subset  \Aa_{I,D}^c$.
\QED

With this lemma the following theorems hold exactly as in \cite{LT}.
\begin{theorem}\label{th:tian1}
The section  $\Phi:B \times \overline{\Ff}_D(M,0,k) \rightarrow \Ee $ defined by equation \ref{phi} gives rise to a generalized Fredholm orbifold bundle with the natural orientation and of index $2c_1(D) +2k+2n-6 +\dim B$.
\end {theorem}
also
\begin{theorem}\label{th:tian2}
Consider two homotopic maps   $i:(B,\p B )\rightarrow (\Aa,\Aa_D^c)$ and  $i':(B',\p B' )\rightarrow (\Aa,\Aa_D^c)$ that represent the same element in $\pi_*(\Aa,\Aa_D^c)$. Suppose also that there  is a deformation between $\omega_B$ and $\omega_{B'}$. If $\Phi'$ and $\Phi$ are defined by equation \ref{phi} for the pairs $(J_{B'},\omega_{B'})$ and $(J_{B},\omega_{B})$ respectively, then $\Phi$ is homotopic to $\Phi'$ as generalized Fredholm orbifold bundles.
\end{theorem}

Theorem \ref{th:tian1} shows that there exist an Euler class 
 $$e([\Phi:B \times \overline{\Ff}_D(M,0,k) \rightarrow \Ee]) \in H_r(B\times \overline{\Ff}_D(M,g,k)$$
 which gives  $\it{the \;\;virtual\;\; moduli\;\; cycle}$
 for the moduli space of equivalence classes of stable parametric holomorphic  maps in the class D,  $\overline \Mm_{0,k}^*(M,D,(J_B,\p J_B))$. This is independent of the choice of $(J_B,\omega_B)$ by \ref{th:tian2}. 

We will denote  by $[\Mm]^{\it{vir}}$ the virtual cycle.
In order to define the invariant we consider the evaluation maps 
$ev_i:B \times \overline{\Ff}_D(M,0,k) \rightarrow M$ given by
$$ev_i(b,[f,\Sigma,x_1,\ldots,x_k])=f(x_i)$$
We then can define
$$PGW^{M,(J_B,J_{\p B})}_{D,0,k}:\bigoplus_{i=1}^k H^{a_i}(M,\Q)^k \rightarrow \Q$$
 by
$$PGW^{M,(J_B,J_{\p B})}_{D,0,k}(\alpha_1,\ldots,\alpha_k)=ev_1^*(\alpha
_1) \wedge \ldots \wedge ev_k^*(\alpha_k)[\Mm]^{\it{vir}}$$
which are zero unless 
\begin{equation}\label{dimmod2}
\sum _{i=1}^k a_i=2c_1(D)+2k+2n-6+\dim B
\end{equation}

We should also point out that if one changes the orientation of $B$ we obtain the same invariant but with  a negative sign.

We have the following theorem:
\begin{theorem}\label{homopgw}
{\bf i)} The  invariants $ PGW^{M,(J_B,J_{\p B})}_{D,0,k}$ are symplectic deformation invariants and  depend only on the relative homotopy class of $(J_B,J_{\p B})$.

{\bf ii)} For  a fixed choice of $k,D$ and $\alpha_i$ the map $\Theta_{0,k,\alpha_1,\ldots,\alpha_k}:\pi_*(\Aa_I,\Aa_{I,D}^c) \rightarrow \Q$, given by\\
$\Theta_{k,\alpha_1,\ldots,\alpha_k}([(J_B,J_{\p B})] = PGW^{M,(J_B,J_{\p B})}_{D,0,k}(\alpha_1,\ldots,\alpha_k)$ is a homomorphism.

\end{theorem}

\proof{} Point $\bf{(i)}$ follows from theorem \ref{th:tian1} and \ref{th:tian2}.
The fact that the morphism $\Theta$ defined in  point $\bf{(ii)}$ is well  defined follows from theorem \ref{th:tian2}. To show that it is a homeomorphism, we choose $(B_1,\p B_1)$, and $(B_2,\p B_2)$ representing 2 elements $\beta_1$ and $\beta_2$ inside $\pi_*(\Aa,\Aa_D^c).$ 
We can choose them such that we can concatenate them through a deformation
 process in which the $J$'s in the intersection do not admit any $\epsilon$ holomorphic   maps in class D. We can therefore see that the new virtual cycle
 corresponding to the classes   $\beta_1+\beta_2 $  will be a disjoint union of the virtual neighborhoods corresponding to $\beta_1$ and $\beta_2$. But this implies that 
the parametric invariants corresponding to the new class $\beta_1+\beta_2 $
 are the sum of the   PGW corresponding to $\beta_1$ and $\beta_2$. Therefore $\Theta$ is a homomorpism.
\QED

We should point out that in the situation that the family $(J_B,J_{\p B})$ satisfies $H_1$ and $H_2$, then the integer valued invariants we obtain by intersecting the image of the pseudocycle \\
$ev_*(\overline {\Mm} _{0,k}^*(M,D,(J_B,\p J_B)))$ with the  classes $(PD(\alpha_1),\ldots,PD(\alpha_k))$  in $H_*(M)^k$ are the same as $PGW^{M,(J_B,J_{\p B})}_{D,0,k}(\alpha_1,\ldots,\alpha_k)$ defined above. Moreover, they can be obtained by counting the number of $J_b$ holomorphic maps in class D with $k$ marked points which intersect  generic cycles representing  $(PD(\alpha_1),\ldots,PD(\alpha_k))$ in $f(z_i)$.

\subsection{A  criterion of parametric regularity for rational maps }\label{subsec22}

  Consider a family of pairs $(J_{b'},\omega_{b'})_{b' \in B'} $ where each $J_{b'}$ is an  almost complex structure tamed by the  symplectic forms $\omega_{b'}$. 

 Since the regularity of a holomorphic map is a local statement within B and it only concerns the almost complex structure data, we claim that for  each $b \in \int B'$ we can restrict our attention to a  neighborhood $B$ such that the family $ J_B$ descends from a fibration in the following sense. We say that the family $(J_B,\omega_B)$ descends from a fibration $M \rightarrow \widetilde{M} \rightarrow B$ if $\widetilde{M}$ comes with 
an almost complex structure  $ \widetilde {J}$ such that the restriction to each fiber $M \times b$ is an almost complex structure $J_b$. Moreover $\widetilde {M}$ admits a closed two form 
$ \widetilde {\omega}$ which restricts on each fiber $M \times b$ to a symplectic form $\omega _b$ that tames $J_b$. Likewise, we can choose a trivialization of the fibration such that smoothly  $\widetilde{M}=B \times M$ and $\pi$ is just the projection on the first factor.
In the following theorem we consider the family of parameters to be a subspace of $\C^m$ and we denote by $z$ the parameters in $\C^m$.
\begin{theorem}\label{th:Mtrans}
  Let  $(J_z,\omega_z)_{z \in \C^m}$  be a family on $M$ descending from the  symplectic
  fibration $(\widetilde{M}, \widetilde{J}, \widetilde{\omega})$
 \begin{equation}
  \xymatrix{
    M \ar[r]^{i} & \widetilde {M} \ar[d]^{\pi}\\
 &\C^m}
  \end{equation}

  Suppose that  $f : \Sigma \longrightarrow M$ is  a $J_0$  holomorphic map and consider the composite map
 $$\widetilde{f}=i \circ f,\widetilde{f} : \Sigma \longrightarrow M \times 0 \subset \widetilde{M}$$ 
 which is $\widetilde{J}$-holomorphic. If $\widetilde{f}$ is regular then   $f$ is $(J_z)$ parametric regular. Moreover, if $\Sigma=S^2$ then the reverse statement holds.
\end{theorem}

The proof of the theorem occupies the rest of the section. Let  $T_{|_{\pi^{-1}(0)}}\widetilde{M}$ be the tangent space along the preimage of $0 \in \C^m$. We will denote by H the subbundle of  $T_{|_{\pi^{-1}(0)}}\widetilde{M}$  which is $\widetilde{\omega}$  orthogonal to the  fiber $\{0\} \times M$. We would like $H$ to  coincide with the  horizontal space of $T\widetilde{M}$ with respect to the trivialization $\pi$ and to be $\widetilde{J}$  invariant. This can be arranged by deforming the form $\widetilde {\omega}$ so that near the zero fiber $\{0\} \times M$ it   is given by 
$$\widetilde{\omega}=\omega_0 + \pi^*(\sigma_{base}),$$
where $\sigma_{base}$  is a standard symplectic two form on the  holomorphic base B. Throughout this deformation process $\widetilde{J}$ is still $\widetilde{\omega}$ tamed.

 Let $g_0$ be a metric on $M_0$ and $\nabla$ be the Levi-Civita connection on $M$ associated
with it. $\nabla^{st}$ will be  the standard
Levi-Civita connection on $\C^m.$ We will denote from now on $\widetilde{\nabla} = \nabla \times \nabla^{st},$
the product connection on $\widetilde{M} \simeq \C^m \times M$. The regularity of $\widetilde{f}:\Sigma \longrightarrow \widetilde{M}$ is by
definition, equivalent to the fact that $D_{\widetilde{f}}$ is surjective,
where $D_{\widetilde{f}}$ is the linearization of $\bar{\partial},$
\begin{equation*}
  \xymatrix{
  D_{\widetilde{f}}:C^{\infty}(\widetilde{f}^* T\widetilde{M}) \ar@{->>}[r]  &
  \Omega_{\widetilde{J}}^{0,1}(\Sigma,\widetilde{f}^* T\widetilde{M}).}
\end{equation*}

Using the connection $\widetilde{\nabla}$ we will derive formulas for
$D_{\widetilde{f}}$ and  express them in terms of the linearization $D \Phi$.

Since $\widetilde{M} \simeq \C^m \times M$ and
$\im \widetilde{f} \subset \{0\} \times M,$ we have the following relations:
$$ \widetilde{f}^* \left(T \widetilde{M} \right)=\widetilde{f}^* \left(T\widetilde{M}_{\pi ^{-1}(0)} \right)=\widetilde{f}^*( H \oplus TM )= triv \oplus f^*(TM ) $$
where by $triv$ we denote the trivial $m$-dimensional complex  bundle
over $\Sigma$. This gives
\begin{equation}\label{sp1}
  C^{\infty}(\widetilde{f}^* T\widetilde{M})  \simeq
  C^{\infty}(triv) \oplus C^{\infty}(f^*TM)
\end{equation}

 Given that each fiber is  $\widetilde{J}$ invariant, and that $H$ is $\widetilde{J}$ invariant along $\pi^{-1}(0)$, we obtain
\begin{equation}\label{sp2}
  \Omega_{\widetilde{J}}^{0,1}(\Sigma,\widetilde{f}^* T\widetilde{M})  \simeq
  \Omega_J^{0,1}\left( \Sigma,f^* TM  \right) \oplus
  \Omega_{\widetilde{J}}^{0,1}(\Sigma,H)
\end{equation}

From \eqref{sp1} and \eqref{sp2} we obtain 
\begin{equation*}
  \xymatrix{
  D_{\widetilde{f}}:C^{\infty}(triv) \oplus C^{\infty}(f^*TM) \ar@{->>}[r]  &
  \Omega_J^{0,1}\left( \Sigma,f^* TM  \right) \oplus
  \Omega_{\widetilde{J}}^{0,1}(\Sigma,H)}
\end{equation*}
and by considering the appropriate restrictions we obtain
the following operators

$$\begin{array}{lccc}
D_{1,vert}: & \Cinf (triv) & \longrightarrow
&  \Omega^{0,1}_{J}\left( \Sigma,f^*TM \right)\\
D_{1,hor}: & \Cinf(triv) & \longrightarrow
   & \Omega^{0,1}_{\widetilde{J}} \left( \Sigma,H \right)\\
D_{2,vert}: & \Cinf \left(f^*TM \right) &  \longrightarrow
  &  \Omega^{0,1}_{J}\left(\Sigma,f^*TM \right)\\
D_{2,hor}: & \Cinf \left(f^*TM \right) & \longrightarrow
 &   \Omega^{0,1}_{\widetilde{J}} \left( \Sigma,H \right)\\
\end{array}$$

We will sometimes use $D_k=(D_{k,vert},D_{k,hor}), \;k=1,2$.

To compute the formulas for these operators we will use the following
general method (see \cite{AA}). Consider $\xi \in C^{\infty}(\Sigma,\widetilde{f}^* T\widetilde{M})$ and
$\widetilde{F}_{\xi}:[0,1] \times \Sigma \longrightarrow \widetilde{M}$ given by
$\widetilde{F}_{\xi}(t,x)  =
  \exp_{\widetilde{f}(x)}^{\widetilde{\nabla}}\left( t\xi(x) \right),$
for $\xi$ sufficiently small. Let $s:\Sigma \longrightarrow T\Sigma$ be a
section and $\widetilde{s}$ its lift to $T\left( [0,1] \times \Sigma \right).$
We denote $\frac{\partial}{\partial t}$ the vector field in
$T\left([0,1] \times \Sigma \right)$ corresponding to the parameter
in $[0,1].$ Define  $\widetilde{f}_t(x) := \widetilde{F}_{\xi}(t,x).$ For any
$x \in \Sigma,$
define the path $\widetilde{\gamma}^{\xi} _x : [0,1] \longrightarrow \widetilde{M}$ given by
$\widetilde{\gamma}^{\xi}_x(t) = \widetilde{F}_{\xi} (t,x),$
the image under $\widetilde{F}_{\xi}$ of $[0,1] \times x$ in $\widetilde{M}.$
By the definition of $\widetilde{F}_{\xi},$
$\widetilde{\gamma}^{\xi}_x$ is a geodesic path in $\widetilde{M}$ relative
to the connection $\widetilde{\nabla}.$ Denote by 
$\tau_{t,x}^{\xi} : T_{\gamma_x(t)}\widetilde{M} \longrightarrow
  T_{\gamma_x(0)}\widetilde{M}$ the parallel transport in $\widetilde{M}$
along the curve $\gamma_x :=\widetilde{\gamma}^{\xi}_x$. To compute $D_{\widetilde{f}}(\xi)(s)$ in
general, one needs to consider the expression
$\frac 12 \tau_{t,x}^{\xi}( d \widetilde{f}_t(s) + \widetilde{J}d \widetilde{f}_t(js))$ and take its derivative with respect to t  at $t=0$ i.e.
\begin{equation}\label{fr:delform}
  D_{\widetilde{f}}(\xi)(s) = \frac 12 \frac{\partial}{\partial t} \left(
  \tau_{t,x}^{\xi}( d\widetilde{f}_t(s) +
  \widetilde{J}d\widetilde{f}_t(js)) \right)_{|_{t=0}}
\end{equation}

We define $Const$ to be the subspace of $C^{\infty}(triv)$ made out  of constant sections.
For the proof of the theorem, we are particularly interested in computing
$D_{1,hor}$ and the restriction of $D_{1,vert}$ to $Const.$ 

 In order to simplify the notation, we denote by $x$ the coordinate
  on $\Sigma$ and write the points  in $\C^m \times M$ as
  $\left( z_1, \ldots, 
    z_m,  y \right) $ where $z_1=w_1+iv_1$ and so on. 
For  simplicity we denote the vector field in $Const$ by  $ \frac{\p}{\p w_k}=\partial_{w_k}$ and so on. Since we are going to work  with an arbitrary choice of $w_k$ and $v_k$  we will refer to them simply as $\partial_{w}$, unless we need  to be more specific.
\begin{lemma}\label{lm:splitexp} The following relations hold: 

  {\bf i)} $D_{2,hor} = 0$ 

  {\bf ii)} $D_{2,vert} = D_f$ 

  {\bf iii)} ${D_{1,hor}}(\xi)=
                 \pbar_{\C^m}(\xi), \forall \xi \in C^{\infty}(triv)$, where $\pbar_{\C^m}$ is the delbar operator in $\C^m.$

  {\bf v)} $({D_{1,vert}})(\p_z)(s) =\frac 12 \frac {\p}{\p z}(J(z))_{|z=0}(df(js)$ for $\p_z$ a typical vector field  in $ Const \subset  C^{\infty}(triv)$.
\end{lemma}

\proof{} Since $\widetilde{f}=f \circ i \subset \{0\} \times M$ we can naturally view any  $\xi \in C^{\infty}(f^* TM)$ as an element in  $C^{\infty}(\widetilde{f}^* T \widetilde{M})$ with values in the  vertical direction tangent to $\{0\} \times M.$ We have that    $\widetilde{F}_{\xi}(t,x) = \exp_{\widetilde{f}(x)}^{\widetilde{\nabla}}
      (t\xi) = \exp_{f(x)}^{\nabla} (t\xi),$ with
    $\im \widetilde{F} \subset \{0\} \times M.$ 
    This implies that the $d\widetilde{f}_t(s)$ are also vertical vector fields
    supported in $\{0\} \times M$ and, since $\widetilde{J}$ keeps
    $T(\{0\} \times M)$ invariant, we have as well that the
    $\widetilde{J} d\widetilde{f}_t(js)$
    are vertical vector fields in $\{0\} \times M.$ 
    Similarly, $\widetilde{F}_{\xi} ^* \frac{\partial}{\partial t}$ is a vertical section in
    $T\widetilde{M}$ supported in $\{0\} \times M$ and parallel transport along
    $\widetilde{f}(x)$ with respect to $\widetilde{\nabla}$ is the same as parallel
    transport with respect to $\nabla.$ 

  A direct application of \eqref{fr:delform} is that
  \begin{equation*}
    \left( D_{\widetilde{f}}\xi \right) (s) = \frac 12
    \frac{\partial}{\partial t}
    \left( \tau_{t,x}^{\xi}d\widetilde{f}_t(s)  +
           \tau_{t,x}^{\xi}\widetilde{J}df_t(js) \right)_{|_{t=0}} 
     = \left(D_f\xi \right)(s),
  \end{equation*}
  which proves ${\bf (i)}$. Relation ${\bf (ii)}$ follows immediately from the formula above, taking into account that $D_{\widetilde{f}} \xi =D_{2,vert}(\xi)$,  and  that     
 $\im D_{\widetilde{f}}|_{C^{\infty}f*TM} \subset  \Omega^{0,1}_{J}\left(\Sigma,f^*TM \right)$. 

  For the proofs of ${\bf (iii)}$ and ${\bf (v)}$ we now consider $\xi \in  C^{\infty}(triv).$ We can assume
  $\xi  = \phi(x)\p _{w}$ where $\phi:\Sigma \leftarrow \C^m$. In this situation,
  $\widetilde{F}_{\xi}(t,x) =
    \exp_{\widetilde{f}(x)}^{\widetilde{\nabla}}(t \p_w) =
    (\phi(x)t,0, \ldots, 0,f(x)).$ It then follows that the paths $\gamma_x$ are
  straight lines in $\C^n \times f( x )\subset  \widetilde{M}$ and therefore
  the parallel transport along ${\gamma}_x,$
  $\tau_{t,x}:T_{(t,f(x))}\widetilde{M} \longrightarrow T_{0,f(x))} \widetilde{M}$
  is the identity.
We are also going to consider the coordinates $x \in \Sigma 
  \text{ of the type } x=x_1 +i x_2$, and do our computations for 
$s=\partial_{x_1}$.

If  $\widetilde{J}(t)$ is the almost complex structure at  $\widetilde{\gamma}^{\xi}_x(t)$ then   
$
\widetilde{J}(t)$ has the form $ \left(
{\begin{matrix}
A_t&0\\
B_t & J_t\\
\end{matrix}}
\right)
$ with respect to the product structure $\C^m \times M.$ Moreover along $\pi^{-1}(0)$   we have 
$
\widetilde{J}(0)=\left( {\begin{matrix}
J_{\C^m}&0\\
0 & J_t\\
\end{matrix}} \right)
$. Therefore  $\frac{\partial}{\partial t}\widetilde{J}(t)$ preserves the fibers, the same as $\widetilde{J}(t)$ does. Moreover, along $\{0\} \times M$, $\widetilde{J}(0)$ preserves the splitting into $TM$ and $H$.
As we have seen parallel transport along $\widetilde{\gamma}^{\xi}_x(t)$ is just the identity.
 
 Considering local coordinates  on $ \Sigma $ $x=x_1 +i x_2$ and taking  $s=\partial_{x_1}$, we have: 
$$\begin{array}{lll}
  D_{1,hor}(\phi \partial _w)(\partial_{x_1}) & = & 
\frac 12  proj_H \frac{\partial}{\partial t}
   \left( \tau_{t,x}^{\xi}d\widetilde{f}_t(\partial_{x_1})  +
 \frac 12  \tau_{t,x}^{\xi}\widetilde{J}d\widetilde{f}_t(j \partial_{x_1})
   \right)_{|_{t=0}} \\ 
  & = & \frac 12 proj_H \frac{\partial}{\partial t}\left( d\widetilde{f}_t
   (\partial_{x_1})  + 
   \frac 12  \widetilde{J}d\widetilde{f}_t( \partial_{x_2}) \right)_{|_{t=0}} \\
  & = & \frac 12 \frac{\partial}{\partial t }
   \left(\partial_{x_1}(\phi (x))t,0, \ldots, 0\right)_{|_{t=0}} +
  \frac 12 proj_H \frac{\partial}{\partial t}
   \left( \widetilde{J}_t\right)_{|{t=0}} df(\partial_{x_2}) + \\
  & & \frac 12 proj_H \widetilde{J}_0 \frac{\partial}{\partial t}
   \left(\partial_{x_2}(\phi(x))  t,0, \ldots,0,df(x)\right)_{|{t=0}} 
\end{array},$$
where, as mentioned before, $\phi:\Sigma \rightarrow \C^m$.
But here the middle term vanishes because  $df(\partial_{x_2})$ is a vertical vector and  $\frac{\p }{\p t} \widetilde J$ preserves fibers so we get that  $ \frac{\partial}{\partial t}
   \left( \widetilde{J}_t\right)_{|{t=0}} df(\partial_{x_2})$ is also a vertical vector. Then
\begin{equation}
 D_{1,hor}(\phi \partial _w)(\partial_{x_1})  = \frac 12 \partial_{x_1} \phi(x) +\frac 12  J_{\C^m}(\partial_{x_2})\phi(x)
\end{equation}

For the last expression we have to use that along $\pi^{-1}(0)$, $\widetilde{J_0}$ preserves the horizontal space $H$, so $proj_H \circ \widetilde{J_0}= \widetilde{J}_{\C^m} \circ proj_H.$
Therefore, the conclusion follows that $D_{1,hor}=\pbar_{\C^m}.$

To prove  point {\bf (v) } of the theorem we  need to consider now $\xi=\partial_w$, that is  $\xi \in Const$. Under this assumption we have  $\tau_{t,x}^{\partial_w}d \widetilde{f}_t=df_0$. Thus  
$$ \frac{\partial}{\partial t} \tau_{t,x}^{\p_w} d\widetilde{f}_t(s) =0.$$
As before, $s$ is a just a section in $T\Sigma$.
We then have
$$\begin{array}{lll}
 D_{1,vert}(\partial _w)(s) & = &\frac 12  proj_V \frac{\partial}{\partial t}\left( \tau_{t,x}^{\p w} d\widetilde{f}_t(s) 
 +\frac 12  \tau_{t,x}^{\p w} \widetilde{J}d\widetilde{f}_t(js) \right)_{|_{t=0}} \\
& = & \frac 12  proj_V\frac{\partial}{\partial t}\left(\tau_{t,x}^{\p w} d\widetilde{f}_t(s)\right)_{|_{t=0}}  
+ \frac12  proj_V\frac{\partial}{\partial t}\left(\tau_{t,x}^{\p w} \widetilde{J}(\tau_{t,x}^{\p w})^{-1} \right)_{|_{t=0}} \cdot df(js) \\
& & +
 \frac 12 proj_V\widetilde{J}_0\left( \frac{\partial}{\partial t} \tau_{t,x}^{\p w} d\widetilde{f}_t(js) \right)_{|_{t=0}} 
 = \frac 12  proj_V\left(\widetilde{\nabla}_{\partial _w} \widetilde{J} \right) df(js) 
\end{array}
$$
where we denote by $proj_V$ the projection onto the fibers. Recall that $\frac{\p}{\p t}\widetilde{J}$ takes vertical vector fields into vertical vector fields.
  Therefore 
 $$\frac 12 proj_V \widetilde{\nabla}_{\partial _w} \widetilde{J} df(js)=
\frac 12 \frac {\p}{\p w}(J(z))(df(js))$$
 precisely because $df(js)$ is a vertical vector field and because the covariant derivative along horizontal vector fields was chosen to be the standard connection in $\C^m$.
Applying the same reasoning for $i \p_v$ we see that 
 $$({D_{1,vert}})(\p_z)(s) =\frac 12 \frac {\p}{\p z}(J(z))_{|z=0}(df(js))$$
 It is worth to point out that  $\frac {\p}{\p z}(J(z)_{|z=0}=d\psi^*_0(\frac {\p}{\p z})$.

 \QED 

\proof{ of the theorem}{\bf Implication ``$\Rightarrow$''}  Using lemma \ref{lm:splitexp}, point {\bf (v)} we get the commutativity of the following diagram

\begin{equation}\label{diag}
\xymatrix{
  T_0\C^m \ar[rr]^{ d \psi} \ar[d]^i & & T_J \Aa_{I} \ar[d]^R \\
  Const \ar[rr]^{D_{1,vert}} & & \Omega_J^{0,1}(\Sigma,f^*TM), }
\end{equation}
where $i:T_0 \C^n \rightarrow Const \subset C^{\infty}(triv)$ is the natural identification map and $\psi$ is the morphism from the parameter space to the space of almost complex structures. $R$ is, as mentioned before, given by $R(Y)=\frac12 Y \circ df \circ j$. 

Since $D_{\widetilde{f}}$ is surjective by hypothesis of, this means that
  $D_1 \oplus D_2$ is  surjective. We therefore have, by lemma (\ref{lm:splitexp}) {\bf (i),(ii)},
  \begin{equation}
  \xymatrix{
    D_1=(D_{1,vert},D_{1,hor}): C^{\infty}(triv) \ar[r] & \coker D_f \oplus
      \Omega_{\widetilde{J}}^{0,1} (\Sigma,H)  }
  \end{equation}
  is surjective. Since  the kernel of the $\pbar_{\C^m}$ operator on $\C^m$ consists precisely of constant sections, lemma \ref{lm:splitexp} {\bf (iii)}    implies
  that $ D_{1,hor}^{-1} (0)=Const.$ Therefore we have that the operator

  \begin{math}
  \xymatrix{
    {(D_{1,vert}})_{|_{Const}} : Const\ar[r] & \coker D_f  }
  \end{math}
  is surjective. But this will imply that 

\begin{math}
  \xymatrix{
    {D_{1,vert}}_{|_{Const}} \circ i :T_0 \C^m \ar@{->>}[r] & \coker D_f  }.
  \end{math}
 But as we saw in the proof of \ref{pr:transreg}, $R$ induces an isomorphism $\widetilde{R}:\widetilde{\coker d\Pi} \longrightarrow 
 \coker D_2$ and moreover the diagram \ref{diag} will be still commutative if we restrict $d\psi$ and $D_{1,vert}$ to $\coker d \Pi$ and $\coker D_2$ respectively.
 Therefore  
 \begin{math}
  \xymatrix{
    d \psi : T_0 \C^n\ar[r] & \coker d \Pi.  }
  \end{math} is surjective.
  By proposition \ref{pr:transreg}, this yields  
 exactly the parametric regularity.

 For the inverse implication, we notice that since $D_{1,hor}$ is $\pbar_{\C^m}$, it will cover the space $ \Omega_{\widetilde{J}}^{0,1} (\Sigma,H) $ when $\Sigma=S^2$. By hypothesis we have that   \begin{math}
  \xymatrix{
    d \psi : T_0 \C^n\ar[r] & \coker d \Pi.  }
  \end{math} is surjective 
and the above observation implies that 
 \begin{equation}
  \xymatrix{
    D_1=(D_{1,vert},D_{1,hor}): C^{\infty}(triv) \ar[r] & \coker D_f \oplus
      \Omega_{\widetilde{J}}^{0,1} (\Sigma,H)  }
  \end{equation} is also surjective. Therefore $D_{\widetilde{f}}$ is a surjective operator.
\QED

\section{Resolutions of singularities and relative PGW}\label{sec3}

\subsection{ Quotient  singularities}\label{subsec31}
In this subsection we will give an overview of work of Kronheimer \cite{K} and Abreu-McDuff \cite{AM} on how to construct special families of almost complex structures arising from the study of the total spaces of deformations for some quotient singularities. In the end of the section we will explain how these families serve our purpose of counting nontrivial PGW.
The local picture is as follows(see Kronheimer \cite{K}):\\
 We consider the  particular type of Hirzebruch-Jung  singularity  $Y_0={\C^2}/{C_{2\ell}}$, given by the   diagonal action by scalars of $C_{2\ell}$ on $\C^2$, where $C_{2\ell}$ is the cyclic group of order $2\ell$. This admits a resolution 
$\sigma_0:\widetilde{Y_0}\rightarrow Y_0$ where $\widetilde{Y_0}$ is the total space of the line bundle of degree $ -2\ell$ over $\C P^1$. 
The exceptional curve of the resolution, we will call it E, is a curve of selfintersection $-2\ell$ and is the zero section of $\widetilde{Y_0}$.
This resolution  admits a $2\ell -1$ complex dimensional  parameter family of deformations
,$\widetilde{Y}_t,t \in \C^{2 \ell-1}$. With the exception of the case $\ell=2$ the total space $\widetilde Y=\bigcup \widetilde Y_t$ of  the family of deformations  is the total space of the  vector bundle ${\cal O}(-1)^{2\ell }$. More precisely, we consider the exact sequence of bundles 
\begin{equation}\label{kroseq}
\begin{array}{ccccccc}
 {\cal O}(-2\ell)&\rightarrow  & {\cal O}(-1)^{2\ell } &\stackrel{r}{\rightarrow}& {\cal O}^{2\ell-1}& \\
\end{array}
\end{equation}
 where $r$ is given by evaluating at $2 \ell -1 $ generic  sections of the dual of $\widetilde{Y}$, $\widetilde{Y}^*={\cal O}(1)^{2 \ell}$. Since holomorphically $ {\cal O } ^{2\ell-1}$ is trivial, we can project it to its fiber $\C^{2 \ell -1}$ and hence we obtain a submersion  $\widetilde{q}: {\cal O}(-1)^{2\ell } \rightarrow \C^{2\ell-1}$ with $\widetilde{Y}_t =\widetilde{q}^{-1}(t)$. Also it can be seen that  $\widetilde{Y}$ is smoothly isomorphic  with $\widetilde{Y_0} \times \C^{2 \ell -1}$ and a choice of trivialization provides an isomorphism 
\begin{equation}
\begin{array}{ccc}
\theta:\widetilde{Y}& \stackrel{C^{\infty}}{ \widetilde{=}}& \widetilde{Y_0} \times \C^{ 2 \ell -1}\\
\end{array}
\end{equation} 

  We consider now a $4 \ell$-dimensional basis of sections in the dual $\widetilde{Y}^*$. Here the space of holomorphic sections is given by
 $ \bigoplus_{i=1}^{2 \ell}H^0(C P^1,{\cal O}(1))\widetilde{=}(\C^2)^{2 \ell}$. Denote by $Y$ the subspace of $(\C^2)^{2 \ell}$ consisting of $2 \ell$ -tuples of vector in $\C^2$ which span either zero or a line. By evaluating all the $4 \ell$ section we obtain a map 
$$ \sigma :\widetilde{Y} \rightarrow Y \subset \C^{4 \ell}$$
 which contracts $E$ to a point $\gamma_0=\sigma(E)$. Moreover,  $\gamma_0$ is the only singular point of $Y$. and the morphism is one to one outside $E$. We  also define a map $q:Y \rightarrow \C^{2 \ell -1}$ by evaluating at the original $2 \ell -1$ generic sections. The following diagram commutes 

\begin{equation}
\xymatrix{
  \widetilde{Y} \ar[r]^{\sigma} \ar[d]^{\widetilde{q}} &  Y \ar[d]^{q}\\
  \C^{2 \ell -1} \ar[r]^{id}  & \C^{2 \ell -1} }
\end{equation}
  We can obtain a 2-form $\tau$ on $Y$ by pulling back a K\"ahler form from  $\C^{4 \ell}$. Via $\sigma^*$ this can be seen as a two form on $\widetilde{Y}$ which restricts to a K\"ahler form $\tau_t$ on each fiber $\widetilde{Y_t}$ if $t \neq 0$ but degenerates along $E$ when $t=0$. If we further push forward through $\theta$, these forms can be seen as a family of forms on the manifold $\widetilde{Y_0}.$

 As in \cite{AM}, we can choose an appropriate compactification of the local picture as follows :

 Let $B^{4 \ell -2}$ be the unit ball in $\C^{2\ell-1}.$  We have a family $(\overline{Y}_t,J_t^{\ell},\tau_t)_{t\in B^{4 \ell -2}}$, where each $(\overline{Y}_t,J_t^{\ell},\tau_t^{\ell}),t \ne 0$ is a K\"ahler manifold diffeomorphic with $S^2 \times S^2$, and, $(\overline{Y}_0,J_0^{\ell})$ is a complex manifold, also diffeomorphic with $S^2 \times S^2$  and $\tau_0$ degenerates along $A-\ell F$. The total space of the family has the following properties:

{\bf a)} The space  $\overline{Y}$
 = $\cup_{t\in B^{4 \ell -2}} \overline{Y}_t$ is smoothly diffeomorphic
 with $S^2 \times S^2 \times B^{4 \ell -2}$. Moreover $\overline{Y}$  is a complex manifold with a complex structure $\widetilde{J^{\ell}}$ which restricts to each fiber $\overline{Y}_t$ to the complex structure $J_t^{\ell}$. Also, $\overline{Y}$  has a closed $(1,1)$ form
$ \tau $ which is satisfies all the properties of a K\"ahler form outside the zero fiber and restricts at each fiber to the forms, $\tau_t$ .

{\bf b)}  The form $\tau$  restricted to $\overline{Y_0}$  degenerates along the exceptional curve $A-\ell F$ 
 
 Since the forms $\tau_t$ are obtained by restricting the closed form  $\tau$  to fibers it is immediate that they are all in the same cohomology class. Moreover, since $(\tau_0)_{| A- \ell F} = 0$ we obtain that $\forall t \in B^{4 \ell -2},[\tau_t^{\ell}]=[\omega_{\ell}]$.

 From {\bf (a)} we see that there is a holomorphic projection $\pi:\overline{Y}\rightarrow S^2 \times B^{4 \ell -2}$. This is because every $\overline{Y}_t$ is a ruled surface therefore it fibers over $S^2$. If we denote be $\alpha$ the area form on $S^2$ we can construct a two form 
$$\tau^{\lambda}=\tau +(\lambda - \ell) \pi^*(\alpha)$$
For $\lambda > \ell$ these forms are K\"ahler forms and moreover they restrict to each $\overline{Y}_t$ to symplectic forms in the class $[\omega_{\lambda}]$.
This proves that any $J_t^{\ell}$ is tamed by a form isotopic with $\omega_{\lambda}$ as long as $\lambda> \ell$.
 We  now follow a similar procedure to construct a family of symplectic forms $\omega_t,t \in B^{4 \ell -2}$ such that each $\omega_t$ tames $J_t^{\ell}$. We will now change the forms $\tau_t$  by perturbing with a a positive factor of  $\pi^*(\alpha)$ only around $t=0$ and smoothen with a cut-off function. By this procedure we obtain symplectic forms $\omega_t$  with variable cohomology classes.

 In conclusion, we have pairs $(S^2 \times S^2,J_t^{\ell},\omega_t)_{t\in B^{4\ell-2}}$ where $\omega_t$ is a {\it symplectic} 
structure  on $S^2 \times S^2$  that tames $J_t^{\ell}$. Moreover $[\omega_t]_{t \in  S^{4\ell-3}}=[\omega_{\ell}]$. This gives a family of almost complex structures which we denote by abuse of notation $B_{\ell}$  such that  $(B_{\ell},\p B_{\ell}) \in (\Aa_{[\ell,\ell+\epsilon]},A_{\ell})$ for any  $\epsilon>0$.
More importantly, only $J_0^{\ell}$ admits the exceptional curve  in the class
$A-\ell F$.

We will then obtain a family  of almost complex structures on 
$(S^2 \times S^2 \times X)$ by taking $(J_t^{\ell}\times J_{st})$, and by 
abuse of notation, we will call this family also $B_{\ell}$. Therefore we 
just produced on $(S^2 \times S^2 \times X)$ pairs $(B_{\ell},\p B_{\ell}) \subset (\Aa_{[\ell,\ell+\epsilon]},A_{\ell})$, with  $ \epsilon>0$ that represents an element $\beta_{\ell}$ in 
$\pi_* (\Aa_{[\ell,\ell+\epsilon]},A_{\ell}) $. Moreover each $B_{\ell} \subset \Aa_{\ell+\epsilon}$ for any small $\epsilon >0$.

From the choice of  the J's we know that the only structure which admits $A-\ell F$ curves is $J_0 \times J_{st}$.

\subsection {The computation of PGW}\label{subsec32}

Here we prove that {\bf $( H_1)$ }  and {\bf $(H_2)$ } are satisfied for the family  $(B_{\ell},\p B_{\ell})$,  and therefore the invariant is integer valued and can be obtained by counting holomorphic maps intersecting generic cycles of appropriate dimension.

{\bf Claim 1. }The family  $(B_{\ell},\p B_{\ell})$ satisfies  {\bf $(H_1)$ }.

\proof{ of claim 1}This is proved by inspection. Only $J_0^{\ell}\times J_{st}$ admits $A- \ell F$ stable  maps, and the only maps in this class are copies of the imbedded map E in any fiber $S^2 \times S^2 \times \point$. Hence there are no decomposable $J_b$ holomorphic maps. We should point out that for other almost complex structures J on $S^2 \times S^2 \times X$ one could have decomposable J-holomorphic maps in the class $A - \ell F$.
\QED

{\bf Claim 2. }  The family  $(B_{\ell},\p B_{\ell})$ satisfies  {\bf $H_2$ }.

\proof{ of claim 2} From the sequence (\ref{kroseq}) we have that the map E, which is $\widetilde J^{\ell}$-holomorphic has the normal bundle  ${\cal O}(-1)^{2\ell }$ and therefore we can apply lemma 3.5.1 pg 38  in \cite{MS} for the  integrable almost complex structure $\widetilde{J}$. If follows that E is $\widetilde J^{\ell}$  regular inside $\overline{Y}.$ If we consider now $\overline{Y} \times X $ and 
$\widetilde J^{\ell} \times J_{st}$, the curve E lies entirely inside  $\overline{Y}$ and therefore the normal bundle inside  $\overline{Y} \times X $ is   ${\cal O}(-1)^{2\ell } \times trivial$, and therefore the curve is is $\widetilde J^{\ell} \times J_{st}$ regular. This splitting and therefore regularity use the fact that the map E is of genus zero. 
Theorem \ref{th:Mtrans}  implies parametric regularity and therefore  {\bf $(H_2)$ } holds.
\QED

We can therefore conclude that the invariants
$$PGW^{S^2 \times S^2 \times X,(B_{\ell},\p B_{\ell})}_{A-\ell F,0,k}:\bigoplus_{i=1}^k H^{a_i}(S^2 \times S^2 \times X,\Q)^k \rightarrow  \Z$$
are integer valued. We have two situations. First, if $X=\point$ then the moduli space of unparametrized curves has dimension 0 so we would count isolated curves. This follows immediately from the fact that $c_1(A- \ell F)=-4 \ell+2$ (adjunction formula) and therefore 
\begin{align*}
  \dim \Mm_{0,0}^*(S^2 \times S^2,A -\ell F,(B_{\ell},\p B_{\ell})) &
  = 2 \times 2  +2c_1(A- \ell F) +\dim B^{\ell} -6 \\
   & = 4 -4 \ell+4  + 4 \ell -2  -6=0.
\end{align*}

Moreover, the invariant $PGW^{S^2 \times S^2 \times X,(B_{\ell},\p B_{\ell})}_{A-\ell F,0,0}([\point])=1$ because it counts E, the only $J_{b_{\ell}}$ map in the class $A- \ell F$.

In the situation that $\dim X=n >0$, we will count maps with one marked point. $c_1(A- \ell F)$ will be the same since the holomorphic maps in class $A- \ell F$ will have the image entirely in the fibers $S^2 \times S^2 \times \point$.
we therefore have 
\begin{align*}
  \dim \Mm_{0,1}^*(S^2 \times S^2 \times X,A -\ell F,(B_{\ell},\p B_{\ell})) &
  = 2 \times ( 2 +n) +2c_1(A- \ell F) +\dim B^{\ell} -6 +2 \\
   & = 2n+2
\end{align*}

 We will consider a cycle  in the homology class $F$ which will lie in a fiber $S^2 \times S^2 \times \point$ inside $S^2 \times S^2 \times X$. It easily follows that the only $J_{b_{\ell}}$ holomorphic map with one marked point which intersect this cycle transversely is a copy of the map E inside the fiber $S^2 \times S^2 \times \point$. We obtain  that
$$PGW^{S^2 \times S^2 \times X,(B_{\ell},\p B_{\ell})}_{A-\ell F,0,1}(PD([F])=1.$$

Applying theorem \ref{homopgw} we obtain that the morphism $\Theta
$ in both situations is nontrivial and therefore there is a nonzero element
 \begin{equation}\label{blcycles}
\beta_{\ell} \in \pi_{4\ell-2}((\Aa_{[\ell,\ell+\epsilon]},A_{\ell}) \text { for all } \epsilon>0
\end{equation}
 that is represented by the cycle $(B_{\ell},\p B_{\ell}) \subset (\Aa_{\ell + \epsilon},\Aa_{\ell+\epsilon,D}^c).$

\section{Almost complex structures and symplectomorphism groups}\label{sec4}
 
\subsection{Almost complex structures and symplectomorphisms; deformations along compact subsets}\label{subsec41}

In this subsection we will give a quick overview of what can be said about the behavior of spaces of almost complex structures and about the symplectomorphisms groups as the symplectic form varies.

We will restrict our attention to variations of the symplectic form $\omega$ along a line L inside $\Kk$, parametrized by the real parameter $\lambda$. 
If L happens to be a ray $\lambda \omega,\lambda >0$ then $G_{\lambda}$ is independent of $\lambda$. It will therefore make sense to consider $L \neq ray$.

If $M=S^2 \times S^2$,  a great deal is known about the structure of  $\Aa_{\lambda}$ see \cite{M3}.
For example, one can establish  that  there is a direct inclusion
$\Aa_{\lambda} \subset \Aa_{\lambda'}$, for $\lambda < \lambda'$. Moreover, the
homotopy type of the spaces $\Aa_{\lambda}$ changes only as $\lambda$ strictly 
passes an integer $\ell$. 

None of this is known to hold for M an arbitrary symplectic manifold. Nevertheless, as a consequence of the fact that taming is an open condition, we are able to establish the following lemma,  which we use in the proof of the theorem \ref{th:Mtrans}
\begin{lemma}\label{lm:kinc}
{\bf i)}  Let $K'$ to be an arbitrary  compact subset of $\Aa_{\lambda}$. Then there is an $\epsilon_{K'}>0$ such that $K'$ is contained in $\Aa_{\lambda+\epsilon}$,  for  $|\epsilon| <\epsilon_{K'}$.

{\bf ii)} 
 Consider $K$ an arbitrary compact set in  $ G_{\lambda}$. Then there is an  $\epsilon_{K}>0$ and a map  $h:[- \epsilon_K,\epsilon _K] \times K \rightarrow \Gg_{|L}$ such that the following diagram commutes

\begin{equation}\label{ext}
\xymatrix{
  h:[- \epsilon_K,\epsilon_K]\times K \ar[r] \ar[d]^{pr_1} &  \Gg_{|L}\ar[d]^{pr_2} \\
  [- \epsilon_K,\epsilon_K] \ar[r]^{incl}  & (- \infty,\infty). }
\end{equation}

Moreover, for any two  such  maps $h$ and $h'$ which coincide on $0 \times K$, there exist, for an $\epsilon' $ small enough, a homotopy $H:[0,1] \times [- \epsilon', \epsilon'] \times K \rightarrow \Gg_{|L}$ between them 
 which satisfies
\begin{equation}\label{germ}
\xymatrix{
 H:[0,1] \times  [- \epsilon ',\epsilon']\times K \ar[r] \ar[d]^{pr_1} &  \Gg_{|L} \ar[d]^{pr_2} \\
 [- \epsilon ',\epsilon'] \ar[r]^{incl}  & (- \infty,\infty). }
\end{equation}

\end{lemma} 
\proof{} Subpoint {\bf (i)} is an immediate consequence of the openness of the taming condition.

For the proof of {\bf (ii)}, let's first notice that, since the symplectic condition is an open condition, there is a convex neighborhood U of $\omega_{\lambda}$ inside the space of 2-forms such that   any closed  $\omega'$  in U is still symplectic. Moreover for any $g_{k_0} \in K \subset G_{\lambda}$
  and any symplectomorphism  $g_k \in K$ which is sufficiently close to $g_{k_0}$ we can choose  $\epsilon(k_0)>0$ small enough, such that $g_k^* \omega_{\lambda + \epsilon}$ is still inside $U,$ for all $0 \leq \epsilon < \epsilon(k_0)$. Since $K$ is compact we can do this process finitely many times such that
 in the end we have an $\epsilon(K) >0$ such that for any $g_k \in K$ $g_k^* \omega_{\lambda + \epsilon} \in \U,\; \text{ forall } \;  0 \leq \epsilon < \epsilon(K) $. We will construct the elements $h(\epsilon,k)$ as follows. For $t \in [0,1]$ the
 forms 
$$\omega_{k,\lambda +\epsilon}^t:=t g_k^* \omega_{\lambda +\epsilon} + (1-t)\omega_{\lambda +\epsilon}$$
 are symplectic since both $g_k^* \omega_{\lambda +\epsilon}$ and $\omega_{\lambda +\epsilon}$ are inside the convex set $U$. 
We now apply Moser's argument  and obtain a family of diffeomorphisms $\xi_{k,\lambda +\epsilon,t}$ with the property that $\xi_{k,\lambda +\epsilon,t}^*\omega_{k,\lambda +\epsilon}^t=\omega_{\lambda +\epsilon}$. We will now
 define $h(\epsilon,k):=g_k \circ \xi_{k,\lambda +\epsilon,1}$. Then $h$ 
has the required properties.

For an arbitrary $h:K \times [-\epsilon, \epsilon] $ satisfying (\ref{ext})  we take a homotopy $F:[0,1] \times  [-\epsilon, \epsilon] \times K \rightarrow \R \times  \Diff M$ given by  $F(t,\epsilon,k):=(\epsilon,h(t \epsilon,k))$. This gives a homotopy between $h$ and  $h_0:[-\epsilon,\epsilon] \times K \rightarrow \R \times \Diff M$, where $h_0(\epsilon',k)=h(0,k)$. We similarly obtain a homotopy $F'$ between $h'$ and $h_0$, where $h'$ also satisfies (\ref{ext}). By concatenating one homotopy  with the opposite of the other we obtain a homotopy between $h$ and $h'$ which we call $G:[0,1] \times  [-\epsilon_1, \epsilon_1] \times K \rightarrow \R \times \Diff M $. Denote by $g_{s,\epsilon,k}:=G(s,\epsilon,k)$. We will now follow the same procedure as before. Namely, we  restrict to a short interval $[ - \epsilon' ,\epsilon']$ such that, if we call
 $$\omega_{s,k,\lambda +\epsilon}^t:=t g_{s,\epsilon,k}^* \omega_{\lambda +\epsilon} + (1-t)\omega_{\lambda +\epsilon}$$ 
 then these are symplectic,  $\forall 0\leq |\epsilon|<\epsilon'$ and $\forall t,s \in [0,1]$. This is possible because $\omega_{s,k,\lambda }^t= \omega_{\lambda}$. Applying Moser's argument again we obtain diffeomorphisms  $\xi_{s,k,\lambda +\epsilon,t}$ with the property that $\xi_{s,k,\lambda +\epsilon,t}^*\omega_{s,k,\lambda +\epsilon,}^t=\omega_{\lambda +\epsilon}$. We will now define $H(s,\epsilon,k):= g_{s,\epsilon,k} \circ \xi_{s,k,\lambda +\epsilon,1}$. Then $H$ has the required properties.   
\QED
\begin{definition}\label{defext}

Let  $\rho:B \rightarrow G_{\lambda}$ be a cycle in  $G_{\lambda}$. An {\bf extension} $\rho^{\epsilon}$ of $\rho$ is a smooth family of cycles $\rho^{\epsilon}:B \rightarrow  G_{\lambda + \epsilon}$ defined for $|\epsilon| \leq \epsilon_0$ such that $\rho^0=\rho$ and satisfying (\ref{germ}). Using \ref{lm:kinc} {\bf (i)} every cycle $\rho$ has an extension.
\end{definition}
\obs\  
  Consider two  extensions $\rho^{\epsilon}_1, 0 \leq |\epsilon| <\epsilon_1$ and  $\rho^{\epsilon}_2, 0 \leq |\epsilon| <\epsilon_2$.  By (\ref{germ}) there is an $\epsilon '>0$  and  a  homotopy  between $\rho^{\epsilon}_1$ and $\rho^{\epsilon}_2$  defined for all $0 \leq \epsilon \leq \epsilon'$. Hence any extension provides  well defined elements in $\pi_*G_{\lambda + \epsilon}$ for small values of $\epsilon$. Therefore each $[\rho] \in \pi_*(G_{\lambda}^X)$ has an  extension $[\rho^{\epsilon}] \in \pi_*(G_{\lambda +\epsilon}^X)$ whose {\bf germ} at $\epsilon =0$ is independent of the choices of $\rho$.

\begin{definition}
We say that a smooth  family of elements $[\rho^{\epsilon}] \in \pi_*G_{\lambda + \epsilon},0<\epsilon < \epsilon_{\rho} $ is {\bf new} if it is not the extension for $\epsilon >0$  of any element $[\rho] \in \pi_* G_{\lambda}$.  
\end{definition}
In the next section we will use the same letter $\rho$ to refer  both to  cycles as well as to the homotopy class they represent.

\subsection{The relation between almost complex structures  and symplectomorphism groups; the role of PGW}

 Consider $(M, \omega_{\lambda})$ symplectic structures on $M$ such that as before, the symplectic forms $\omega_{\lambda}$ span a line in  a positive cone ${\cal K}$ inside $H^2(M,\R)$. Denote by 
\begin{equation}\label{ellplus}
\Aa_{\ell^+}=\{ J \mid \text{ there is an } \epsilon_{J}>0 \text{ s.t. } J \in \Aa_{\ell+\epsilon} \text{ forall } 0<\epsilon<\epsilon_J \}
\end{equation}
 
 \begin{definition}
 Consider a nontrivial element  $\beta_{\ell} \in \pi_* (\Aa_{\ell^+},\Aa_{\ell})$. We say that $\beta_{\ell}$ is a {\bf persistent} element if its image under the natural morphism

$$i_* \pi_*(\Aa_{\ell^+},\Aa_{\ell}) \rightarrow  \pi_*(\Aa_{[\ell,\ell + \epsilon]},\Aa_{\ell})$$ is nonzero for any $\epsilon$ arbitrary small.
 \end{definition}
\proof{ of the theorem \ref{th:Mts}} 
We will  consider the long exact sequence of relative homotopy groups of the pair $(\Aa_{\ell^+},\Aa_{\ell})$ 

 \begin{math}
  \xymatrix{
    \ldots \ar[r] & \pi_{k} \Aa_{\ell^+} \ar[r] &  \pi_{k}( \Aa_{\ell^+},\Aa_{\ell})  \ar[r] & \pi_{k-1} \Aa_{\ell} \ar[r] &  \pi_{k-1} \Aa_{\ell^+} \ar[r] & \ldots   }
  \end{math}

Since by construction $\beta_{\ell} \in \pi_{k}( \Aa_{\ell^+},\Aa_{\ell}) $ is nontrivial, then one of the two 
following cases can happen:\\
{\bf Case 1} $\beta_{\ell} \mapsto  \gamma_{\ell} \ne 0 \in \pi_{k-1} \Aa_{\ell}$\\
{\bf Case 2}  $\beta_{\ell} \mapsto 0 \in \pi_{k-1} \Aa_{\ell}$. In this situation, there is an element $ 0 \ne\alpha_{\ell}  \in \pi_{k} \Aa_{\ell^+}$ where $\alpha_{\ell} \mapsto \beta_{\ell}.$

We will do the analysis case by case for our situation:\\
{\bf Case 1 } If we are in this case then  we consider the fibration
  (\ref{fr:Dseq}), that yields

\begin{math}
\xymatrix{
 & G_{\ell} \ar[r] & \Diff_0(M) \ar[r] & \Aa_{\ell}}
\end{math}

We consider the long exact sequence in homotopy
                    
\begin{math}
  \hspace{-0.2in}
  \xymatrix
  {
    &\ldots \ar[r] & \pi_{k-1}(G_{\ell} ) \ar[r] &
      \pi_{k-1}\Diff_0(M) \ar[r] & &  \\
    &\ar[r] & \pi_{k-1} \Aa_{\ell} \ar[r] & \pi_{k-2} G_{\ell}  \ar[r] &
     \pi_{k-2}\Diff_0(M) \ar[r] & \ldots 
  }
\end{math}

Again, there are two possibilities:

{\bf i) }$\gamma_{\ell} \rightarrow \theta_{\ell} \ne 0 \in \pi_{k-2}G_{\ell} .$
In this situation, we have a nontrivial element 
 $ \theta_{\ell} \in \pi_{k-2}G_{\ell}$, such that
 $\theta_{\ell} \mapsto 0  \in \pi_{k-2}\Diff_0(M)$.  Then we are in case {\bf A}.\\
 This element is {\it fragile}. This can be proved by contrapositive. Assume that  $\theta_{\ell}$ can be extended by $\theta_{\ell + \epsilon}$ which yields nontrivial classes in $\pi_{k-2}G_{\ell + \epsilon}$. Then  $\theta_{\ell+\epsilon} \mapsto 0  \in \pi_{k-2}\Diff_0(M)$ as well. Therefore it appears as a boundary of an  element $\gamma_{\ell+\epsilon} \in \pi_{k-1}\Aa_{\ell + \epsilon}$ which is homotopic with $\gamma_{\ell}$. But by construction and lemma (\ref{lm:kinc}), we know that $\gamma_{\ell}$ is a contractible cycle inside $\Aa_{\ell + \epsilon}$. This contradicts the existence of $\gamma_{\ell + \epsilon}$.\

\obs\ {\bf  Direct construction of the elements $\eta_{\ell+ \epsilon}$ }  \label{conscyc}

 If the fragile element $\theta_{\ell} \pi_{k-2}G_{\ell}$ is trivial, then the elements $\eta_{\ell+ \epsilon} \in \pi_{k-1}G_{\ell +\epsilon}$ are constructed as follows.

 Since $\theta_{\ell}$ is trivial, it can be represented by a cycle $\theta_{\ell}$ which is the boundary of a $k-1$-dimensional disc $D_{\ell} \in G_{\ell}$. By lemma (\ref{lm:kinc}) point{\bf (ii)} we can extend $D_{\ell}$ to discs  $D_{\ell+\epsilon}$ inside $\G_{\ell+\epsilon}.$ 
Alternatively, we can push by extensions the cycles $\theta_{\ell}$ into cycles $\theta_{\ell +\epsilon}$ inside $G_{\ell+\epsilon}$. 
These  cycles are null homotopic therefore bound discs $C_{\ell+\epsilon}$ inside $G_{\ell+\epsilon}$. Due to (\ref{germ}), one can see that $\theta_{\ell+\epsilon}$ and $\p D_{\ell+\epsilon}$ are homotopic for small $\epsilon$, therefore we can glue $C_{\ell+\epsilon}$ and $D_{\ell+\epsilon}$ along their boundaries, and obtain a cycle which we denote by $\eta_{\ell+\epsilon}$. 
In what will follow we basically show that if $\theta_{\ell}$ is trivial then $\eta_{\ell+\epsilon}$ gives nontrivial elements in homotopy.
\QED\\
{\bf ii) } $\gamma_{\ell} \mapsto  0 \in \pi_{k-2}G_{\ell}$. Then
$\gamma_{\ell}$ is in the image of the morphism 
 $ \pi_{k-1}\Diff_0(M )\rightarrow 
   \pi_{k-1} \Aa_{\ell}$, and therefore there is an element,
 $ \gamma_{\ell}' \in  \pi_{k-1}\Diff_0(M )$ such that 
$0 \ne \gamma_{\ell}' \mapsto \gamma_{\ell}$.

In this situation, we can choose a cycle $S \subset \Aa_{\ell}$ representing
$\gamma_{\ell} \in \pi_{k -1}(\Aa_{\ell})$, and, using lemma (\ref{lm:kinc}), there is an $\epsilon_S >0$ 
such that for any $\epsilon \text { such that } 0<\epsilon < \epsilon_S,  
  S \subset \Aa_{\ell +\epsilon}$. We make the following\\
{\bf Claim (1)} $0=[S] \in \pi_{k -1} \Aa_{\ell+\epsilon}$.

 By hypothesis $S$ is the boundary of a cycle $B_{\ell}$ such that $B_{\ell} \subset \Aa_{\ell+\epsilon}$ for all small $\epsilon>0$. Therefore we have a $k$ dimensional ball inside $\Aa_{\ell+\epsilon }$ whose boundary is $S$, which proves the claim.
We therefore have:

\begin{math}
  \hspace{2in} \gamma'_{\ell} \longrightarrow \;\;\; [S]=0 \in \pi_{k -1} \Aa_{\ell+\epsilon }
\end{math}

\begin{math}
  \hspace{-0.4in}
  \xymatrix
  {
    &\ldots \ar[r] & \pi_{k-1}(G_{\ell+\epsilon})  \ar[r] &
       \pi_{k-1}\Diff_0(M) \ar[r] \ar@^{=}[d] &
       \pi_{k-1} \Aa_{\ell+\epsilon} \ar[r] &
       \pi_{k-2}(G_{\ell+\epsilon} ) \ar[r] & \ldots  \\
    &\ldots \ar[r] & \pi_{k-1}(G_{\ell} ) \ar[r] &
       \pi_{k-1}\Diff_0(M) \ar[r] &
       \pi_{k-1} \Aa_{\ell} \ar[r] \ar@^{.>}[u]_{i_{|k}}&
       \pi_{k-2}(G_{\ell} ) \ar[r] & \ldots 
  }
\end{math}

\begin{math}
  \hspace{2in} \gamma'_{\ell} \longrightarrow \;\;\; \gamma_{\ell} \in \pi_{k -1} \Aa_{\ell+\epsilon }
\end{math}

Here, from the first row, since $\gamma'_{\ell}$ is in the kernel of the map $ \pi_{4\ell-3}\Diff_0(M) \rightarrow 
  \pi_{k-1} \Aa_{\ell+\epsilon}$, it has to be in the image of the map
$ \pi_{k-1}(G_{\ell+\epsilon} ) \rightarrow  \pi_{k-1}\Diff_0(M)$, and therefore we are able to produce an element 
$0 \ne \eta_{\ell + \epsilon} \in  \pi_{k-1}(G_{\ell+\epsilon} )$ such that $\eta_{\ell +\epsilon}$  persists in the topology of the group of diffeomorphisms.
Thus we are in case {\bf B}.

 The elements we obtain here are {\it new}. This follows easily by assuming the opposite. That is, if we consider that there is an element $0 \neq \eta_{\ell} \in \pi_{k-1}G_{\ell}$ whose germ is given by $\eta_{\ell + \epsilon}$, then the image of $\eta_{\ell}$ in $\Diff_0(M)$ has to be $\gamma'_{\ell}$. But this contradicts the fact that  $\gamma'_{\ell} \mapsto \gamma_{\ell} \neq 0$.\\
{\bf Case 2}. In this situation we have a nontrivial element 
$\alpha_{\ell} \in \pi_{k}\Aa_{\ell^+}$. We then have the following :\\
{\bf Claim (2) } There is an $\epsilon$ such that for $0<\delta < \epsilon $  $\alpha_{\ell}$ has a representative $C$ inside  $ \Aa_{\ell + \delta}$, $0 \ne [C] \in \pi_{k}\Aa_{\ell + \delta}$.
The proof of this statement follows from the construction of $\alpha_{\ell}$.
Namely, since $\beta_{\ell} \mapsto 0 \in \pi_{k-1}\Aa_{\ell}$ we conclude that there exist a $k$-dimensional disk D inside $\Aa_{\ell}$ whose boundary is $\p B_{\ell}$; by lemma (\ref{lm:kinc}) {\bf (i)} this can be viewed inside $\Aa_{\ell + \delta}$ for small $\delta$. We can now glue $B_{\ell}$ and $D$ along their boundary $\p B_{\ell}$. In this manner we get a cycle $C \subset \Aa_{\ell+\delta}$ which represents the class $\alpha_{\ell}$.
We can therefore consider again the sequence

\begin{math} 
  \hspace{-0.2in}
  \xymatrix
  {
    \ldots \ar[r] & \pi_{k}(G_{\ell +\delta}) \ar[r] &
      \pi_{k}\Diff_0(M) \ar[r] & &  \\
    \ar[r] & \pi_{k} \Aa_{\ell+\delta} \ar[r] & \pi_{k-1} G_{\ell+\delta} \ar[r] &
     \pi_{k-1}\Diff_0(M) \ar[r] & \ldots 
  }
\end{math}

\noindent{\bf Claim (3)}   $[C]$ doesn't lift to a nontrivial 
element in $\pi_{k}\Diff_0(M) $.

\proof{}{of claim (3)} We should first make the observation that there is a map
\begin{equation}  \pi_{k}\Diff_0(M) \rightarrow \pi_{k} \Aa_{\lambda}
\end{equation}
for any $\lambda$ and moreover as $\lambda$ varies this maps are homotopic in $\Aa_I$. If $C$  did lift, the map $\pi_{4\ell -2}\Diff_0(M) \rightarrow 
 \pi_{k}\Aa_{\ell}$ would produce a cycle $[B]$ $ \in \Aa_{\ell}$,
which by means of lemma (\ref{lm:kinc}) can be viewed inside all 
$\Aa_{\ell+\epsilon}$ for small $\epsilon$  and which moreover is  homotopic with C inside  $ \Aa_{[\ell, \ell +\epsilon ]}$. Therefore $[C]$ would map to $0 \in \pi_{k}(\Aa_{\ell^+},\Aa_{\ell})$, which would contradict its definition.
\QED\\
Since $[C]$ cannot be in the image of the map $\pi_{k}\Diff_0(M) \rightarrow \pi_{k}\Aa_{\ell+\delta}$, we know that [C] must have nonzero image   $[C] \mapsto \eta_{\ell + \delta} \ne 0$ in $\pi_{k-1}G_{\ell+\delta}$. Moreover form the obvious properties of exact sequences again, $ \eta_{\ell + \delta} \rightarrow 0$ through the natural inclusion map $\pi_{k-1}G_{\ell+\delta} \rightarrow \pi_{k-1}\Diff_0(M)  $. The fact that this elements are {\it new} follows again by assuming the opposite. If they would form  the  {\bf germ} of an element  $\eta_{\ell}$ in $\pi_{k-1}G_{\ell}$, then  $\eta_{\ell}$  would also be null homotopic inside $\Diff_0(M)$ so it would therefore  come from a class $[C']$ in $ \pi_{k}\Aa_{\ell}$. Moreover, $C'$ would be homotopic with $  C$ inside $\Aa_{[\ell,  \ell + \delta]}$ therefore also in $(\Aa_{[\ell, \ell + \delta]},  \Aa_{\ell})$ which is false given that $C$ has to yield a nontrivial element in $\pi_{ k}(\Aa_{[\ell, \ell + \delta ]},  \Aa_{\ell})$.  Thus we are in the  case {\bf B} of the theorem.

With this, we have exhausted all the possible cases given by the nontrivial PGW.       
\QED\\

 Assume that there is an $\ell$ such that there is no $J$ in $\Aa_{\ell}$ which can be represented by a $J$-holomorphic curve in the class $D$. 
 Then we have the following proposition
\begin{proposition}\label{pers}
Assume that no $J$ in $\Aa_{\ell} $ admits J- holomorphic stable maps in class D. Consider an element $0 \ne  \beta_{\ell} \in \pi_* (\Aa_{\ell^+},\Aa_{\ell})$ obtained by counting nontrivial parametric Gromov-Witten invariants. Then $\beta_{\ell}$ is a persistent element.
\end{proposition}

The proof follows directly from  the theorem (\ref{homopgw}).
\QED\\
Now consider the manifold $(S^2 \times S^2 \times X,\omega_{\lambda} \oplus \omega_{st})$. As explained in (\ref{blcycles})  the cycles $(B^{\ell},\p B^{\ell})$ satisfy the definition (\ref{ellplus}), so they give by Prop (\ref{pers}) {\bf persistent} elements in $\pi_{4 \ell -2}(\Aa_{\ell^+},\Aa_{\ell})$. 

 Therefore theorem (\ref{th:Mts}) applies and so the corollary (\ref{cr:Mts}) holds.

\end{document}